\documentclass[12pt]{article}

\usepackage{amsmath,amssymb,cite,color,url}
\usepackage{fullpage}
\usepackage[pdftex]{graphicx}
\usepackage{multirow}
\usepackage{tikz}

\newcommand{\BEAS}{\begin{eqnarray*}}
\newcommand{\EEAS}{\end{eqnarray*}}
\newcommand{\BEA}{\begin{eqnarray}}
\newcommand{\EEA}{\end{eqnarray}}
\newcommand{\BEQ}{\begin{equation}}
\newcommand{\EEQ}{\end{equation}}
\newcommand{\BIT}{\begin{itemize}}
\newcommand{\EIT}{\end{itemize}}
\newcommand{\BNUM}{\begin{enumerate}}
\newcommand{\ENUM}{\end{enumerate}}

\newcommand{\eg}{{\it e.g.}}
\newcommand{\ie}{{\it i.e.}}

\newcommand{\reals}{{\mbox{\bf R}}}

\makeatletter
\long\def\@makecaption#1#2{
   \vskip 9pt
   \begin{small}
   \setbox\@tempboxa\hbox{{\bf #1:} #2}
   \ifdim \wd\@tempboxa > 5.5in
        \begin{center}
        \begin{minipage}[t]{5.5in}
        \addtolength{\baselineskip}{-0.95pt}
        {\bf #1:} #2 \par
        \addtolength{\baselineskip}{0.95pt}
        \end{minipage}
        \end{center}
   \else
    \hbox to\hsize{\hfil\box\@tempboxa\hfil}
   \fi
   \end{small}\par
}
\makeatother

\newcounter{oursection}

\newcounter{lecture}

\newcounter{algorithmctr}[section]
\renewcommand{\thealgorithmctr}{\thesection.\arabic{algorithmctr}}
{\refstepcounter{algorithmctr}\begin{list}{}{%
\setlength{\rightmargin}{0\linewidth}%
\setlength{\leftmargin}{.05\linewidth}}%
\rmfamily\small
\item[]{\setlength{\parskip}{0ex}\hrulefill\par%
\nopagebreak{\bfseries\textsf{Algorithm \thealgorithmctr~}}}}%
{{\setlength{\parskip}{-1ex}\nopagebreak\par\hrulefill} \end{list}}

\newlength\LabelWidth

\begin{document}

\title{Incremental Proximal Multi-Forecast\\ Model Predictive Control}

\author{Xinyue Shen \and Stephen Boyd}

\maketitle

\begin{abstract}
Multi-forecast model predictive control (MF-MPC) is a 
control policy that creates
a plan of actions over a horizon for each of a given set of 
forecasted scenarios or contingencies,
with the constraint that the first action in all plans be the same.
In this note we show how these coupled plans
can be found by solving a sequence of single plans, using
an incremental proximal method.
We refer to this policy as incremental proximal model predictive 
control (IP-MPC).
We have observed that even when the iterations in IP-MPC are not 
carried out to convergence, we obtain a policy that
achieves much of the improvement of MF-MPC over single-forecast model predictive
control (MPC).
\end{abstract}

\section{Multi-forecast model predictive control}

We consider a control policy which generates a control
input $u_t \in \reals^m$ from the current system state $x_t\in \reals^n$
(presumed known), as well as other information that is available 
in (discrete time) period $t$.

\subsection{Model predictive control}
In model predictive control (MPC), we first form an approximation of 
the dynamics of the system over from period $\tau = t$ to 
period $ \tau= t+H$, where $H$ is the planning horizon.
The approximate dynamics are affine, of the form
\[
x_{\tau+1} = A_{\tau|t} x_{\tau} + B_{\tau|t} u_\tau + c_{\tau|t},
\quad \tau = t, \ldots, t+H-1.
\]
Here $x_t$ is the current state, which is known; $x_{t+1}, \ldots, x_{t+H}$ is our
plan for the future states. The current control input (that we seek) is $u_t$;
$u_{t+1}, \ldots, u_{t+H-1}$ is our plan for future control actions.
The data that define the dynamics used in our plan,
$A_{\tau|t}$, $B_{\tau|t}$ and $c_{\tau|t}$, are a forecast or prediction,
based on information known at period $t$.
(These can change with $t$, as new information becomes available.)
We let $x=(x_{t+1}, \ldots, x_{t+H})$ and 
$u=(u_t, \ldots, u_{t+H-1})$ denote the state and control action plans, 
respectively.

In MPC we choose the plan by solving the optimization problem
\BEQ\label{e-mpc}
\begin{array}{ll}
\mbox{minimize} & G_t(x_t,x,u)\\
\mbox{subject to} & 
x_{\tau+1} = A_{\tau|t} x_\tau + B_{\tau|t} u_\tau + c_{\tau|t}, 
\quad \tau = t, \ldots, t+H-1,
\end{array}
\EEQ
with variables $x \in \reals^{Hn}$ and $u\in \reals^{Hm}$,
where $G_t:\reals^n \times \reals^{Hn} \times \reals^{Hm}
\to \reals \cup \{\infty\}$ is a convex cost function.

Infinite values of $G_t$ are used to enforce constraints, such as 
$u_\tau \in \mathcal U_\tau$, a convex set of allowed control actions, 
or $x_{t+H} \in \mathcal X^\text{term}$, a convex set of allowed terminal states
for our plan.
Traditionally the cost function $G_t$ is separable across $(x_\tau,u_\tau)$,
but it need not be; for example it could penalize a maximum excursion of the 
state or maximum actuator use over the horizon.
Like the dynamics, the cost function $G_t$ can also depend on forecasts of unknown 
future quantities, based on information available at period $t$.
The model predictive control problem \eqref{e-mpc} is convex, 
and readily solved \cite{BoydVand04}, even in real time 
\cite{wang2009fast,mattingley2012cvxgen,stellato2020osqp,banjac2017embedded,Domahidi2013ecos}.
The total number of scalar variables in the MPC problem \eqref{e-mpc}
is $H(m+n)$. 

\paragraph{History.}
MPC has a long history and large literature, and is widely used.
Some early work is \cite{richalet1978model,cutler1980dynamic,garcia1989model}; 
for more recent surveys see the papers 
\cite{morari1999model,holkar2010overview,lee2011model,yu2013model,mayne2014model,abughalieh2019survey} 
or books \cite{maciejowski2002predictive,grune2011nonlinear,camacho2013model,rawlings2017model,rakovic2018handbook}.
Papers describing applications of MPC in specific areas include
data center cooling \cite{lazic2018data},
building HVAC control \cite{afram2014theory,drgovna2020all},
wind power systems \cite{hovgaard2015model},
microgrids \cite{hu2021},
pandemic management \cite{CARLI2020373,peni2020nonlinear},
dynamic hedging \cite{primbs2009dynamic},
railway systems \cite{felez2019model},
aerospace systems \cite{eren2017model},
and agriculture \cite{ding2018104}.
With appropriate forecasting (which in many applications is very simple) and 
choice of cost function, MPC can work well, even though it does not explicitly
take into account uncertainty in the dynamics and cost, or more precisely, 
since it is based on a single forecast of these quantities.
(It does have recourse, however, since the forecasts and plans are updated in
each time period.)

\subsection{Multi-forecast model predictive control}
There are many extensions of MPC that attempt to improve performance 
by taking into account uncertainty in the future dynamics and cost.
Examples include robust MPC \cite{campo1987robust,bemporad1999robust},
min–max MPC \cite{raimondo2009min},
tube MPC \cite{mayne2005robust},
and stochastic MPC \cite{mesbah2016stochastic,heirung2018stochastic}.

One particularly simple approach is multi-forecast MPC
(MF-MPC) which replaces the single forecast of $G_t$, $A_{\tau|t}$, $B_{\tau|t}$,
and $c_{\tau|t}$ 
used in MPC with multiple forecasts, all of which are considered plausible.
We denote these as 
\[
G^i_t, ~A_{\tau|t}^i, ~B_{\tau|t}^i, ~ c_{\tau|t}^i, 
\quad \tau = t, \ldots, t+H-1, \quad 
i = 1, \ldots, S,
\]
where the superscript $i$ gives the forecast or scenario,
and $S$ is the number of scenarios.
We can also specify positive weights $w^1, \ldots, w^S$
associated with these forecasts,
which are often taken to be one.
These multiple forecasts can be found several ways.  Each could be a 
forecast using a different but reasonable method; or they can be samples from
a statistical model of future values. In the latter case, options include
Monte Carlo sampling \cite{shapiro2003monte}, 
Monte Carlo sampling with importance sampling \cite{glynn1989importance},
pseudo-Monte Carlo sampling \cite{caflisch1998monte},
or sigma points in unscented transform \cite[\S 19]{kochenderfer2022algorithms}.
The scenarios can also be constructed by hand.
It is common to take scenario~$i=1$ as the
single forecast that would be used in basic MPC.

To find the desired control action $u_t$, we solve the following 
extension of the single-forecast MPC problem \eqref{e-mpc}:
\BEQ\label{e-mf-mpc}
\begin{array}{ll}
\mbox{minimize} & \sum_{i=1}^S w^i G_t^i(x_t,x^i,u^i)\\
\mbox{subject to} & 
x_{\tau+1}^i = A^i_{\tau|t} x_\tau^i + B_{\tau|t}^i u_\tau^i + c_{\tau|t}^i, 
\quad \tau = t, \ldots, t+H-1, \quad i=1, \ldots, S\\
& u^i_t = u_t, \quad i=1, \ldots, S,
\end{array}
\EEQ
with variables $x^i \in \reals^{Hn}$, $u^i\in \reals^{Hm}$, and 
$u_t\in \reals^m$ (our control action).
In MF-MPC we create $S$ different future state and action plans,
each using its own forecast of cost and dynamics, and add the constraint
that all plans must agree on the first action $u_t$.
Like the MPC problem \eqref{e-mpc}, the MF-MPC problem \eqref{e-mf-mpc} 
is also convex, but it can be a large problem if $S$ is large,
with a total number of scalar variables $SH(m+n)+m$.

\paragraph{History and related literature.}
MF-MPC is a simple special case of many other methods
for incorporating uncertainty and information patterns via multiple scenarios
\cite{golub1995stochastic,zenios1998dynamic,
de2005stochastic,topaloglou2008dynamic,
lucia2012new,maiworm2015scenario,
rockafellar2019progressive,de2021risk}.
In these papers (and others) a scenario tree is used to 
represent the evolution of uncertainty over time,
with non-anticipativity constraints imposed so
that inputs from the same tree node are equal.
Such multi-stage problems can then be solved by methods such as
nested Benders decomposition \cite{birge1988multicut,de2005stochastic},
progressive hedging \cite{rockafellar1991scenarios}, 
progressive decoupling \cite{rockafellar2019progressive},
scenario decomposition with alternating projections \cite{de2021risk},
to evaluate the policy.
In this context, MF-MPC is the very special case where there 
the scenario tree consists of the root 
(the current period), and $S$ edges to the different scenarios at period $t+1$.
In terms of stochastic control, the information pattern for which 
MF-MPC is optimal is one where there are only $S$ possible outcomes 
(the scenarios): When the first action is taken, the $S$ outcomes are 
known but which one will obtain or realize is not;
after the first action is taken, which of the scenarios is realized
is revealed.

The MF-MPC problem is a basic and standard two-stage stochastic programming
problem; see, \eg, 
\cite{sand2004modeling,engell2009online}
\cite[\S 5]{moehle2019dynamic}.
In this context the first stage is referred
to as the ``here and now'' decision and the
second stage as the ``recourse'' actions.
For linear two-stage programs
there are several algorithms that decompose the problem
into each scenario based on cutting plane techniques
\cite{van1969shaped,ruszczynski1986regularized,
sen1994solution,sen1994network,prekopa1995two,
linderoth2001implementing,ruszczynski2003decomposition},
and for convex (nonlinear) two-stage problem there
are methods based on augmented Lagrangians and ADMM
\cite{rockafellar1991scenarios,arpon2020admm}  
\cite[\S 7.6]{parikh2014proximal}.

In summary, neither MF-MPC nor special methods for 
solving the MF-MPC problem, which is a two-stage stochastic 
programming problem, are new. 

\subsection{This note}
The point of this note is to describe a method for solving the MF-MPC
problem \eqref{e-mpc} using an incremental proximal method,
a sequential algorithm that in each iteration
solves a problem similar to the single-forecast MPC problem \eqref{e-mpc}.
If these iterations are continued long enough, the method will converge
to a solution of the MF-MPC problem.   Of more practical interest,
we have found that stopping the incremental proximal algorithm early,
well before it has converged, we obtain a policy that works
well in practice, yielding most or all of the benefits of MF-MPC over
single-forecast MPC.

In this note we are not concerned with comparing MF-MPC with 
other control policies, or arguing that it is a good policy.
Our only point is that MF-MPC can be 
(approximately) evaluated using an iterative method in which each iteration
is essentially solving a (single forecast) MPC problem.
Thus the MF-MPC policy can be evaluated by solving a sequence of MPC problems,
each one asssociated with one scenario.

\section{Incremental proximal model predictive control}
\subsection{Incremental proximal method}
We define $F^i:\reals^m \to \reals \cup\{\infty\}$,
$i=1, \ldots, S$, as the optimal value of the MPC problem in scenario $i$,
as a function of the first action $u_t$.
Specifically, $F^i(u_t)$ is the optimal value of the problem
\[
\begin{array}{ll}
\mbox{minimize} & w^i G_t^i(x_t,x^i,u^i)\\
\mbox{subject to} & 
x_{\tau+1}^i = A^i_\tau x_\tau^i + B_\tau^i u_\tau^i + c_\tau^i, 
\quad \tau = t, \ldots, t+H-1\\
& u^i_t = u_t,
\end{array}
\]
with variables $x^i$ and $u^i$.
This is a convex function since it is the partial
minimization of a convex function over some variables 
(here, $x^i$ and $u^i_{t+1}, \ldots, u^i_{t+H-1}$)
\cite[\S 3.2.5]{BoydVand04}.
In terms of $F^i$, the MF-MPC problem \eqref{e-mf-mpc} can be expressed as the problem
\BEQ\label{e-mf-mpc-Fi}
\begin{array}{ll}
\mbox{minimize} & \sum_{i=1}^S F^i(u_t),
\end{array}
\EEQ
with variable $u_t \in \reals^m$.
(This is the same as the MF-MPC problem \eqref{e-mf-mpc}, after we optimize 
over the variables $x^i$ and $u^i_{t+1}, \ldots u^i_{t+H-1}$.)

The incremental proximal method 
\cite{bertsekas2011incremental,bertsekas2011survey}
solves problems with the sum form in \eqref{e-mf-mpc-Fi}.
In the $k$th iteration the updated iterate $u_t^{(k+1)}$ 
is the solution of the problem
\BEQ\label{e-incre-prox-update}
\begin{array}{ll}
\mbox{minimize} & \alpha_k F^{i_k}(u_t) + \frac{1}{2} \|u_t - u_t^{(k)}\|_2^2,
\end{array}
\EEQ
where $u_t$ is the variable, and $\alpha_k > 0$ is a step size.
The sample index $i_k$ can be chosen in a cyclic order
($i_k = k \mod S$),
or uniformly randomly drawn from the $S$ scenarios.
The step sizes should be square summable but not summable, with a typical 
choice $\alpha_k = \alpha/(k+\beta)$, where $\alpha$ and $\beta$
are positive parameters.
The update \eqref{e-incre-prox-update} is the proximal operator
of $\alpha_k F^{i_k}$ 
\cite{rockafellar1976monotone,Lemaire1989,parikh2014proximal},
giving the method its name.
(`Incremental' refers to the fact that the scenarios are handled separately, one
in each iteration.)

Solving the problem~\eqref{e-incre-prox-update} is the same as solving
the problem
\BEQ\label{e-incre-prox-update-unfold}
\begin{array}{ll}
\mbox{minimize} & \alpha_k w^{i_k} G_t^{i_k}(x_t, x, u) + \frac{1}{2} 
\|u_t - u_t^{(k)}\|_2^2,\\
\mbox{subject to} & 
x_{\tau+1} = A_{\tau|t}^{i_k} x_\tau + B_{\tau|t}^{i_k} u_\tau 
+ c_{\tau|t}^{i_k}, \quad \tau = t, \ldots, t+H-1,
\end{array}
\EEQ
where $u$ and $x$ are the variables.
Solving this convex optimization problem gives $u_t^{(k+1)}$.
The problem \eqref{e-incre-prox-update-unfold} is identical to 
the single forecast MPC problem, except for the addition of the 
proximal term $\frac{1}{2}\|u_t-u_t^{(k)}\|_2^2$.
In many cases the cost function $G_t$ already includes a quadratic
term in $u_t$, in which case
the problem \eqref{e-incre-prox-update-unfold} has exactly the same
form as the single forecast MPC problem.
Thus the cost of carrying out 
each iteration of the incremental proximal method is essentially
the same as the cost of evaluating a single forecast MPC policy.
We refer to the policy obtained by running a fixed number of iterations
of the update \eqref{e-incre-prox-update-unfold} as the incremental
proximal model predictive control (IP-MPC) method.

\paragraph{Convergence.}
While the convergence proof given in
\cite[\S 3]{bertsekas2011incremental} does not exactly cover 
our case here, it is readily modified to handle it, with a few additional
assumptions.
Suppose $G^i$ has the form
\BEQ\label{e-cvg-condition-G}
G^i(x_t, x^i, u^i) = g^i(x_t, x^i, u^i) + \mathbb{I}(u^i \in \mathcal{U}^i)
+ \mathbb{I}(x^i \in \mathcal{X}^i),
\EEQ
where $g^i$ is a real-valued convex function,
$\mathbb{I}$ is the indicator function,
the set $\mathcal{X}^i$ is convex,
and the set $\mathcal{U}^i$ is convex and bounded.
In this case, with step sizes that are square summable but not 
summable, we can conclude that $u_t^{(k)}$ converges to a 
solution of the MF-MPC problem \eqref{e-mf-mpc}.
To see this, observe that when $G^i$
satisfies \eqref{e-cvg-condition-G},
the corresponding $F^i$ is a sum of
a real-valued function and the indicator function
of a bounded set.  Though the bounded set
depends on $i$, while the constraint set
in \cite{bertsekas2011incremental} does not,
with a minor modification one can show that
the same convergence conclusions hold.

While it is nice to know that if the iterations of IP-MPC were continued 
indefinitely, we would (asymptotically) solve the MF-MPC problem \eqref{e-mf-mpc},
an essential part of IP-MPC is that performs well (as a policy) even when 
it is terminated long before it has
solved the MF-MPC problem \eqref{e-mf-mpc} to high accuracy.

\subsection{Mini-batch IP-MPC}

A useful extension of IP-MPC uses multiple scenarios in each update, 
\ie, a minibatch of scenarios, rather than just one.
Thus in each iteration we solve not a signle-forecast MPC problem but a smaller
MF-MPC problem, with $b \ll S$ scenarios, where $b$ is the minibatch size.

In the $k$th iteration, we take a subset of indices
$\mathcal{S}_k \subset \{1, \ldots, S\}$ with $|\mathcal{S}_k|=b$,
which can be chosen cyclically or drawn randomly
from the scenario indices, with or without replacement.  Then the updated
iterate $u_t^{(k+1)}$ is the solution of the problem
\[
\begin{array}{ll}
\mbox{minimize} & 
(\alpha_k/b) \sum_{i \in \mathcal{S}_k}
F^{i}(u_t) + \frac{1}{2} \|u_t - u_t^{(k)}\|_2^2
\end{array}
\]
with variable $u_t$, which is an MF-MPC problem with $b$ scenarios, plus
the quadratic proximal term.
This problem can be expressed as
\[
\begin{array}{ll}
\mbox{minimize} & 
(\alpha_k/b) \sum_{i \in \mathcal{S}_k}
w^{i} G_t^{i}(x_t, x^i, u^i) + \frac{1}{2} 
\|u_t - u_t^{(k)}\|_2^2,\\
\mbox{subject to} & 
x_{\tau+1}^i = A_{\tau|t}^{i} x_\tau^i + B_{\tau|t}^{i} u_\tau^i 
+ c_{\tau|t}^{i}, \quad \tau = t, \ldots, t+H-1, \quad i \in \mathcal{S}_k \\
& u_t = u_t^i, \quad i \in \mathcal{S}_k,
\end{array}
\]
where the variables are $x^i$, $u^i$,
for $i \in \mathcal{S}_k$ and $u_t$.
The number of scalar variables is
$bH(n + m) + m$, which is about $b$ times of that of
the IP-MPC problem \eqref{e-incre-prox-update-unfold}, and a factor $S/b$ smaller
than the full MF-MPC problem \eqref{e-mf-mpc}.

As the minibatch size $b$ increases, the
incremental proximal method converges faster (in terms of number of iterations),
but the iterations are more costly.  A good choice of $b$ trades off these two
competing trends.

\section{Example}

In this section we illustrate IP-MPC with a simple energy storage arbitrage problem.

\subsection{Problem and policies}

\paragraph{The problem.}
We are to choose the charging (discharging, when negative) rate
of an energy storage system in each hour, with time-varying 
energy prices, so as to maximize our average profit.
We let $u_t \in \reals$ denote the battery charging rate in period (hour) $t$,
and $q_t$ the stored energy.  These must satisfy $-D \leq u_t \leq C$ 
and $0 \leq q_t \leq Q$ for all $t$, where $D$ is the maximum discharge 
rate, $C$ is the maximum charge rate, and $Q$ is the storage capacity.
The storage dynamics is given by $q_{t+1}=q_t + u_t$.
The cost in period $t$ is given by
$p_t(u_t + \eta |u_t|)$, where $p_t \geq 0 $ is the mid-price and 
$\eta \in (0,1)$ gives a gap between the buy and sell prices, \ie,
we purchase energy at higher price $(1+\eta)p_t$, and we sell it back at
the lower price $(1-\eta)p_t$.
The charging rate $u_t$ is chosen with knowledge of the current stored 
energy $q_t$ and current price $p_t$, but not future prices $p_{t+1}, \ldots$,
which, however, can be forecast.
The goal is to minimize the average cost, \ie, to maximize the 
average profit.

In this problem the state is the stored energy $q_t$, and the actual system
dynamics are linear, constant, and known.  The true cost is convex, and 
in this case separable across periods.
The only uncertainty is in the future energy prices, which affects the cost.

\paragraph{MPC policy.}
We use a planning horizon $H=24$, \ie, one day.
We denote the forecast of future prices as $\hat p_\tau$, $\tau=t+1, 
\ldots, t+23$.  For notational simplicity, we use $\hat p_t = p_t$, 
the known current price.
At time $t$, we plan the input $u_\tau$ from $\tau=t$ to $\tau=t+23$, \ie,
we plan over the next 24 hours.
To determine $u_t$ we solve the problem
\[
\begin{array}{ll}
\mbox{minimize} & \sum_{\tau=t}^{t+23} \left( \hat p_\tau u_\tau + 
\eta \hat p_\tau |u_\tau| \right)\\
\mbox{subject to} &
-C \leq u_\tau \leq D, \quad \tau = t, \ldots, t+23\\
& 0  \leq q_\tau \leq Q, \quad \tau = t+1, \ldots, t+24\\
& q_{\tau+1} = q_\tau + u_\tau, \quad \tau = t, \ldots, t+23\\
& q_{t+24} = Q / 2,
\end{array}
\]
with variables $u_t,\ldots, u_{t+23}$ and $q_{t+1}, \ldots, q_{t+24}$.
(The current stored energy $q_t$ is known.)
The terminal constraint requires that in our plan, the terminal storage energy 
should be half the capacity.

\paragraph{MF-MPC policy.}
We denote the $S$ forecasts of future prices as
$\hat p_\tau^i$, $\tau=t+1, \ldots, t+23$, $i=1,\ldots, S$,
with $\hat p_\tau^1$ the forecast used in MPC.
As in MPC, we take $\hat p_t^i = p_t$, the known current price.
We create plans $u_\tau^i$, $\tau=t, \ldots, t+23$,
$i=1\ldots, S$, with the constraint that $u_t^1 = \cdots = u_t^S$,
with the common value giving us $u_t$.
We solve the problem
\[
\begin{array}{ll}
\mbox{minimize} & (1/S) \sum_{i=1}^S 
\sum_{\tau=t}^{t+23} \left( \hat p_\tau^i u_\tau^i
+ \eta \hat p_\tau^i|u_\tau^i| \right)\\
\mbox{subject to} &
-C \leq u_\tau^i \leq D, \quad \tau = t, \ldots, t+23, \quad i=1, \ldots, S\\
& 0  \leq q_\tau^i \leq Q, \quad \tau = t + 1, \ldots, t+24, \quad i=1, \ldots, S\\
& q_{\tau+1}^i = q_\tau^i + u_\tau^i, 
\quad \tau = t, \ldots, t+23, \quad i=1, \ldots, S\\
& q_{t+24}^i = Q / 2,\quad i=1, \ldots, S\\
& u_t^1 = \cdots = u_t^S,
\end{array}
\]
with variables $u_\tau^i$, $\tau=t, \ldots, t+23$, $i=1, \ldots, S$, and
$q_\tau^i$, $\tau=t+1, \ldots, t+24$, $i=1, \ldots, S$.
We take $q_t^i = q_t$, the current known stored energy.

\subsection{Parameters and data}

\paragraph{Parameters.}
We take $C=D=10$ and $Q=50$, so we can completely charge or discharge 
our storage system in $5$ hours.  We take $\eta=0.075$, which means 
there is a 15\% difference between the energy buy and sell prices.

\paragraph{Price data.}
We use real price data,
the hourly verified real-time local marginal price (LMP), in dollars per MWh,
for zonal node 51217 obtained from the PJM market \cite{PJM},
over a period of $268$ weeks from July 2016 through August 2021.
We clip or winsorize the smallest values at the 0.2-percentile of prices,
which is $6.6$; the maximum price over this time period was $690$.
The mean price is $28.5$, and the median price is $23.9$.

\begin{figure}
\centering
\includegraphics[width=\textwidth]{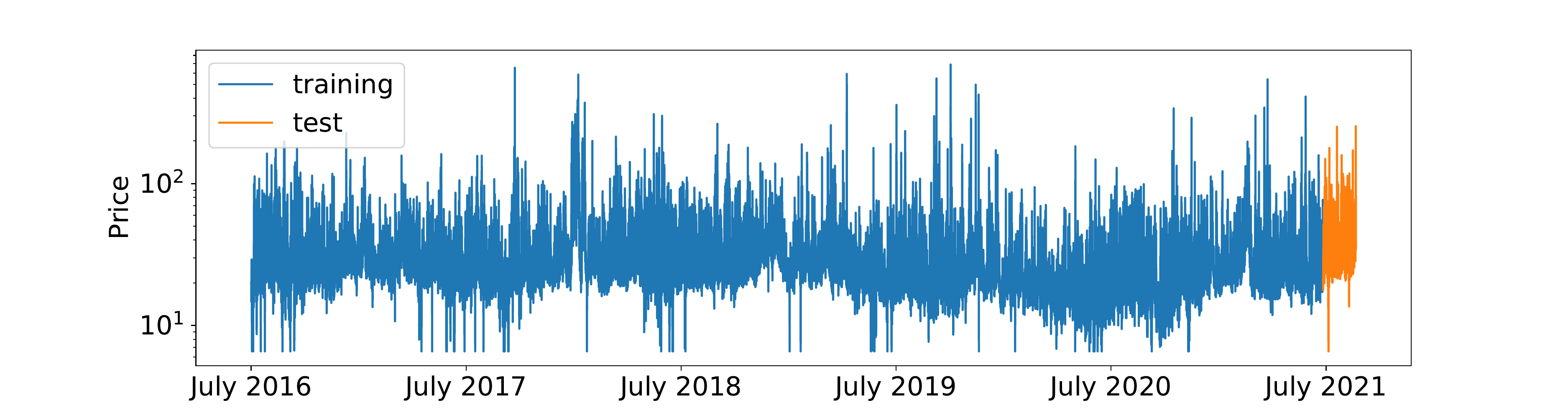}
\caption{Electricity price $p_t$, in dollars per MWh, over $268$ weeks. 
Prices shown in blue are used to fit the forecasting models; the prices
shown in orange are used to evaluate various policies.} 
\label{fig:wholesale_price}
\end{figure}

We use the data of the first $260$ weeks,
from July 1 2016 to June 24 2021, to fit our forecasting model,
and the data of the last $8$ weeks,
June 25 to August 19 2021, to evaluate our policies.
The prices are shown in figure~\ref{fig:wholesale_price}, with blue
showing the data used to develop our forecast model, and orange showing 
the price data used to evaluate the policies.

\subsection{Forecasts}
\paragraph{Transformation.}
The price data are very right skewed, so we first transform them with
\emph{two} log transforms, and work with $z_t = \log \log p_t$, which results in
a distribution of values that is reasonably Gaussian.
While our forecasting methods will use $z_t$, we convert our 
predictions (single or multiple) back to
prices using $\hat p_t = \exp \exp \hat z_t$.

\paragraph{Baseline model.}
We use the simple forecast method described in 
\cite[\S A]{moehle2019dynamic}, applied to $z_t$, which first fits
a baseline to the price data that captures the daily, weekly, and seasonal
variation.
Our baseline $b_t$ has the form
\[
b_t = \beta_0 + \sum_{i=1}^{16}
\left(\beta_i \cos(2\pi t / T_i)
+ \alpha_i \sin(2\pi t / T_i) \right),
\]
where the model parameters are
$\beta_0$, $\beta_i$ and $\alpha_i$
for $i=1, \ldots, 16$.
The periods $T_i$ are as follows.
\begin{itemize}
\item \emph{Diurnal (daily) variation}: $T_i = 24 / k$, $k=1,2,3,4$
\item \emph{Weekly variation}: $T_i = 7 \times 24 / k$, $k=1,2,3,4$
\item \emph{Seasonal (annual) variation}: $T_i = 365 \times 24 / k$, $k=1,2,3,4$
\item \emph{Interaction terms}:
$T_i = 7 \times 24 \pm 24$, $365 \times 24  \pm 24$
\end{itemize}
Thus our basic daily, weekly, and seasonal variation models each have $4$ 
Fourier coefficients; the interaction terms allow the baseline daily and weekly
patterns to vary (a bit) over the year.

\begin{figure}
\centering
\includegraphics[width=\textwidth]{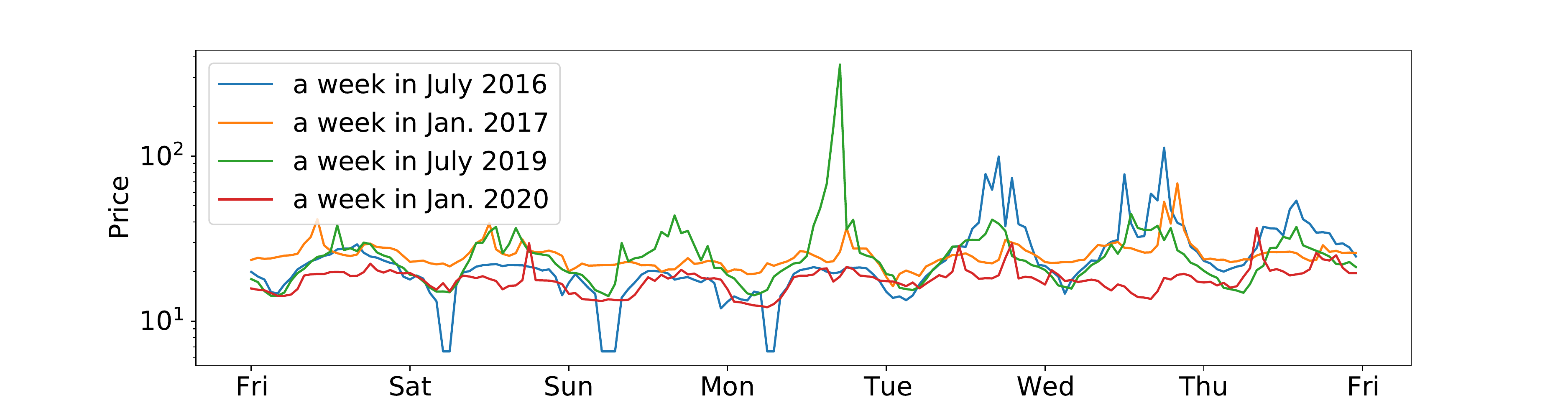}
\includegraphics[width=\textwidth]{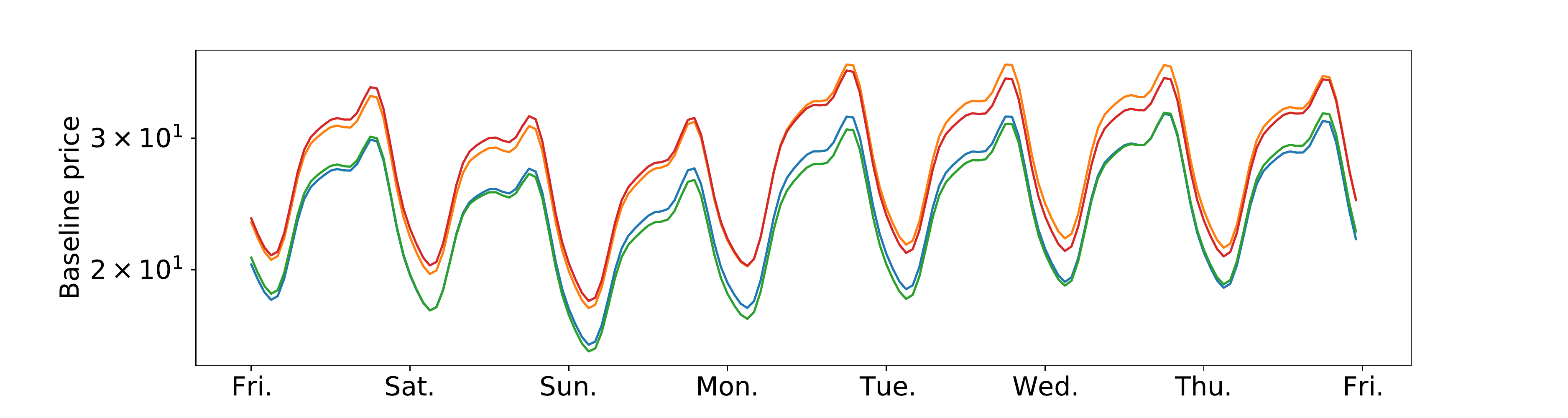}
\caption{Price $p_t$ and baseline price $\exp \exp b_t$,
over four different weeks,
in July 2016, January 2017, July 2019, and January 2020.}
\label{fig:example_week}
\end{figure}

Our baseline model has $33$ parameters; we fit these on the $43680$ data points
using ridge regression.  The log-price prediction error 
$\log p_t - \exp b_t$ has RMS value $0.38$, meaning the baseline
typically differs from the actual price 
by a factor of $\exp 0.38 = 1.46$, \ie, $46\%$. (Recall that the original
price data varies over a 100:1 range.)
The price $p_t$ and baseline price $\exp\exp b_t$
over four different weeks are shown in figure~\ref{fig:example_week}.
Comparing the vertical scales of the two plots,
we can see that the baseline does not capture the occasional large deviations, 
low or high, in the actual prices.

\paragraph{Single forecast model.}
Following the simple forecasting method from 
\cite[\S A]{moehle2019dynamic}, we fit an
auto-regressive (AR) model to predict the
residual $r_t = \log\log p_t - b_t$ over the next $23$ time periods,
given the previous $24$.
The residual AR model has the form
\[
(\hat{r}_{t+1|t},\ldots,\hat{r}_{t+23|t})
= \Gamma (r_{t-23},\ldots,r_{t}),
\]
where $\Gamma \in \reals^{23 \times 24}$ is the AR parameter matrix.
We fit $\Gamma$ using ridge regression on $43633$ training data points.
Our final (single) price forecast is then given by
\[
\hat p_{\tau|t} = \exp \exp \left(b_\tau + \Gamma_{\tau - t} (r_{t-23},\ldots,r_{t}) \right),
\]
where $\Gamma_{\tau - t}$ is the $(\tau-t)$th row of $\Gamma$.
These forecasts have an RMS log-price error of $0.30$, a reduction from
the baseline RMS log-price error $0.38$.
This means our forecasts are typically off from the true price by
around $35\%$.
(If our forecasts were much better, then there would no need to use MF-MPC instead
of MPC.)

\begin{figure}
\centering
\includegraphics[width=\textwidth]{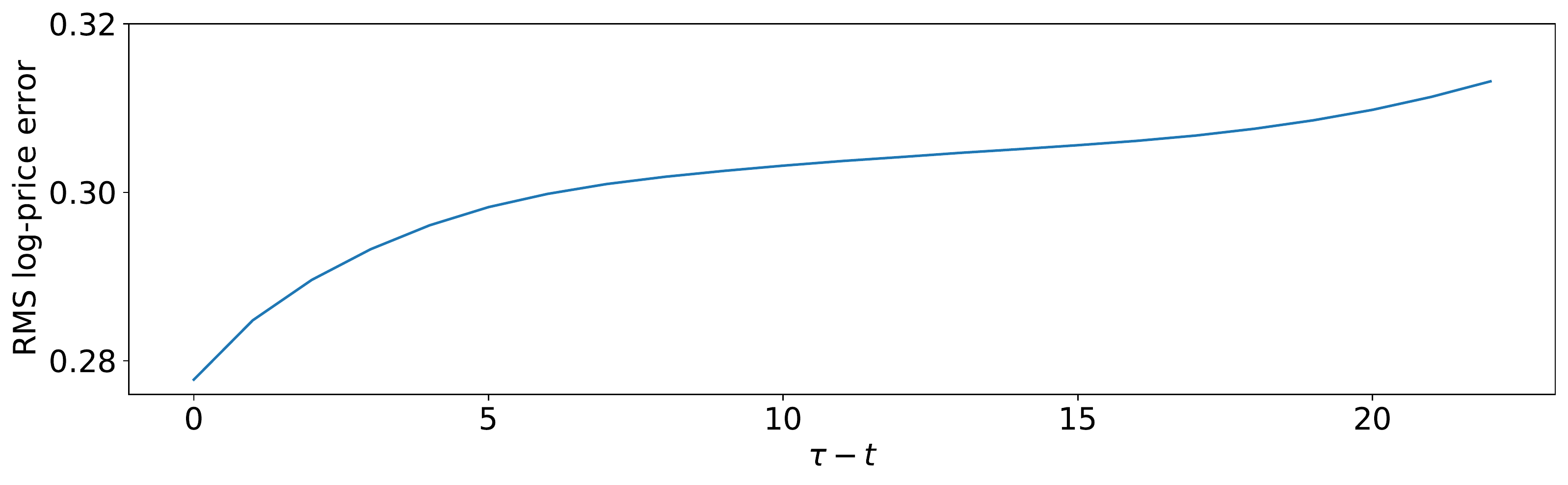}
\caption{RMS forecast error for log-price versus $\tau-t$.}
\label{fig:forecast_error}
\end{figure}

The forecast error varies with $\tau-t$, the number of hours forward that we 
are predicting.
The RMS log-price error versus $\tau-t$ is shown in figure~\ref{fig:forecast_error}.
We can see that our prediction of the next hour's price (\ie, $\tau-t=1$) 
is typically around $32\%$. For larger prediction horizons it increases.
Our predictions 24 hours in the future (\ie, $\tau-t=23$) have
RMS log-error $37\%$, still well below the RMS log-error of the 
baseline alone, which is $46\%$.

\paragraph{Multi-forecast model.}
Denote the AR forecast error at $\tau$
given start time $t$ as 
\[
e_{\tau|t} = \log\log p_{\tau} -\log\log\hat{p}_{\tau|t}.
\]
The statistics of these errors varies with the start time hour of the week,
so we fit a different mean and
a different covariance matrix for each hour of the week.
The $24 \times 7 = 168$ mean vectors 
denoted as $\mu_m$, $m=1,\ldots, 168$, are empirical means.
We fit $168$ covariance matrices to the errors, 
denoted $\Sigma_m$, $m=1,\ldots, 168$, using a Laplacian regularized stratified
model \cite{tuck2021fitting,tuck2021eigen}.
Our graph on the stratified variable, 
in this case hour of the week, is a cycle graph with $24 \times 7$ vertices.
We use the solver implemented in
\cite{tuck2021distributed} to fit the models.

\begin{figure}
\centering
\includegraphics[width=\textwidth]{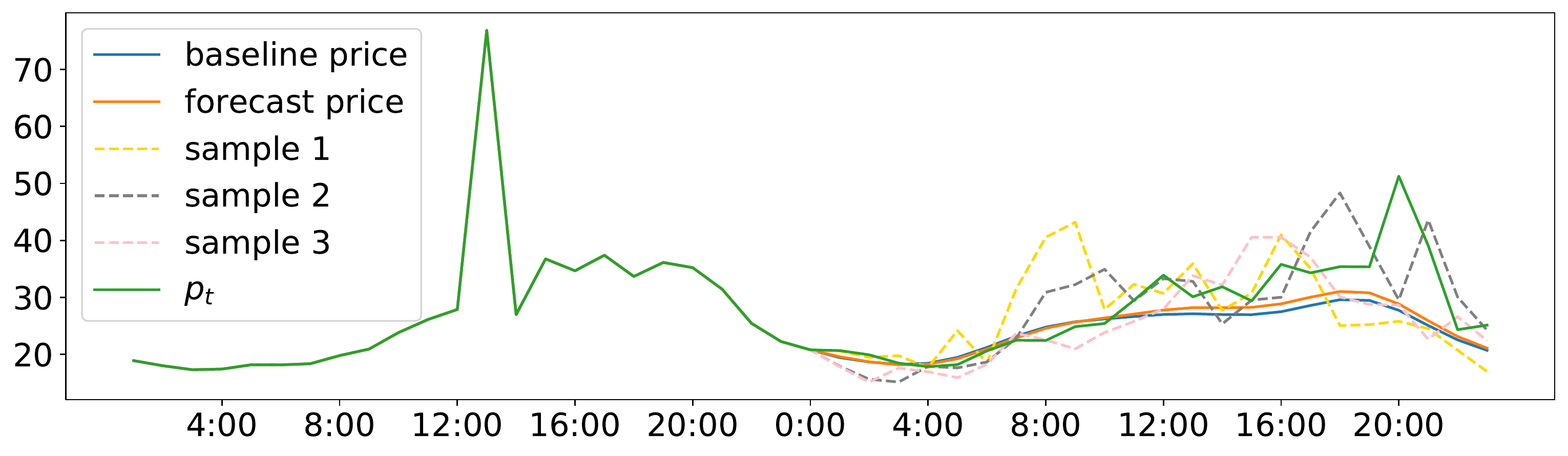}
\caption{Baseline price, forecast price, 
and three price samples in the first 24 test hours,
and the price data from the last 24th training hour 
to the first 24th test hour.
}
\label{fig:forecast_sample_example}
\end{figure}

\paragraph{Scenario sampling.}
Given a start time $t$, $S$ samples 
$(e_{t+1|t}^i, \ldots, e_{t+23|t}^i)$, $i=1,\ldots,S$
are generated IID from $\mathcal{N}(\mu_m, \Sigma_m)$, 
where $m$ is the hour of the week for $t$.
From these we obtain our price forecasts as
\[
\hat{p}_{\tau|t}^i = \exp\exp\left(b_\tau +
\Gamma_{\tau - t} (r_{t-23},\ldots,r_{t}) + e_{\tau|t}^i\right),
\quad
\tau=t+1,\ldots,t+23, \quad i=1,\ldots,S.
\]
Figure \ref{fig:forecast_sample_example} 
shows the actual price, baseline, single forecast, and three random sample forecasts
in the first 24 test hours.

\subsection{Simulation results}

\paragraph{Policies.}
We simulate a number of different policies over the 8 week period at 
the end of our data.
\begin{itemize}
\item MPC using the single forecast.
\item MF-MPC with $S=20, 40, 80, 160$, $320$, and $640$.
\item IP-MPC, with minibatch size $b=20$, for 
a number of iterations $1, 2, 4, 8, 16, 32$.
\end{itemize}
The MF-MPC policies (with different numbers of scenarios)
and IP-MPC use the same samples in the same order, so with $1$ iteration,
IP-MPC uses the same set of samples as MF-MPC with $S=20$ samples, and so on.
For IP-MPC, we start at the MPC plan, so with zero iterations,
this coincides with simple MPC.
We use step sizes $\alpha_k = 7/k$.

\paragraph{Prescient bound.}
To get a performance bound, we compute the exact optimal charging
with all future prices known. (This is a single LP that extends over
the test period.)
The resulting cost is $-62.3$ per hour, \ie, we make an averge profit
of $\$ 62.30$ per hour.
This is an upper bound on how well any policy can do.

\begin{figure}
\centering
\includegraphics[width=\textwidth]{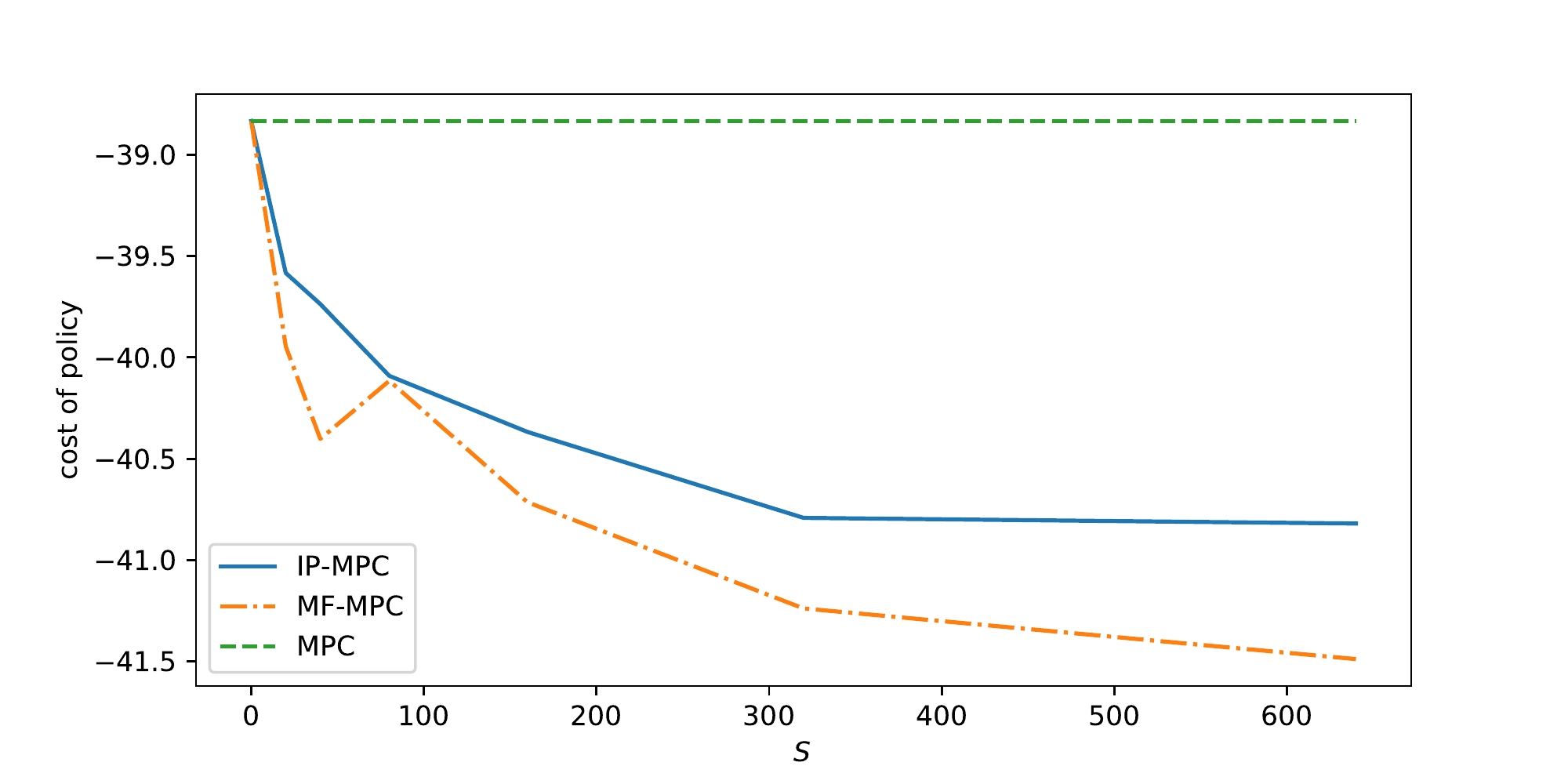}
\caption{
Cost per hour of various policies, including MPC, 
MF-MPC with different numbers of scenarios $S$, and IP-MPC for 
different numbers of iterations. 
The per hour cost of the prescient policy with all future 
prices known is $-62.3$.}
\label{fig:cost_policy_iters}
\end{figure}

\paragraph{Results.}
We ran MF-MPC and IP-MPC for $10$ trials,
with different sets of randomly generated scenarios,
and we report the averaged result over the $10$ trials.
The cost per hour for our policies is
shown in figure~\ref{fig:cost_policy_iters}.
We can see that MF-MPC gives an improvement over MPC,
with $S=640$ samples reducing cost from around $-38.8$ to around $-41.5$,
around $2/3$ of the optimal average cost from our prescient bound.
We can also see that MF-MPC with $S=640$ scenarios is not significantly better
than $S=320$, and that $S=160$ gives us a reasonable fraction of
the improvement over MPC.

The figure also shows that IP-MPC performs well.
For $4$ iterations, IP-MPC uses the first $80$ samples and achieves
a cost that is not far from MF-MPC $S=80$, despite $4$ iterations of IP-MPC
giving only a very crude approximate solution of the associated MF-MPC problem.
For $16$ iterations IP-MPC uses the first $320$ samples and achieves
performance not too far from MF-MPC with $S=320$ samples.
Running IP-MPC for even $32$ iterations does not solve the batch MF-MPC planning
problem to high (or even modest) accuracy, but it yields a policy that does very well.

\paragraph{Computation times.}
We use CVXPY \cite{cvxpy_paper,agrawal2018rewriting} 
and the solver ECOS \cite{Domahidi2013ecos}.
Simulating the MPC policy
over 8 weeks (\ie, $1344$ hours) takes $33$ seconds, which means
around $0.025$ seconds per policy evaluation.
Simulating MF-MPC with 640 samples takes 2268 seconds,
corresponding to $1.69$ seconds per policy evaluation.
Using disciplined parameterized programming \cite{agrawal2019differentiable},
running 32 IP-MPC iterations takes 1023 seconds,
corresponding to $0.76$ seconds per policy evaluation.

\section{Conclusion}

IP-MPC is an iterative method for evaluating an MF-MPC policy, 
with each iteration involving the solution of a single forecast MPC problem.  
In the limit as the number of iterations increases, IP-MPC coincides with MF-MPC.
More interesting to us, and evident in our example, is that IP-MPC
can deliver much of the benefit of MF-MPC with a modest number of 
iterations, well before the IP-MPC policy is close to MF-MPC.

Our example suggests a reasonable general design approach.  First, design a 
standard (single forecast) MPC control policy, choosing the objective and
constraints, and tuning their parameters, to achieve good closed-loop
performance, ideally on real data that was not used to develop the forecaster.
Many variations on the objective and constraints can be tried out,
since the policy is relatively fast to evaluate.
Then, try out MF-MPC with varying numbers of scenarios to see what 
improvement (if any) is obtained over MPC, while possibly making 
small changes to the parameters to improve performance.
Finally, try out IP-MPC, experimenting with the step length sequence
and number of scenarios.

\clearpage

\bibliography{refs}

\begin{thebibliography}{10}

\bibitem{abughalieh2019survey}
K.~Abughalieh and S.~Alawneh.
\newblock A survey of parallel implementations for model predictive control.
\newblock {\em IEEE Access}, 7:34348--34360, 2019.

\bibitem{afram2014theory}
A.~Afram and F.~Janabi-Sharifi.
\newblock Theory and applications of {HVAC} control systems--a review of model
  predictive control ({MPC}).
\newblock {\em Building and Environment}, 72:343--355, 2014.

\bibitem{agrawal2019differentiable}
A.~Agrawal, B.~Amos, S.~Barratt, S.~Boyd, S.~Diamond, and J.~Kolter.
\newblock Differentiable convex optimization layers.
\newblock In {\em Advances in Neural Information Processing Systems}, pages
  9558--9570, 2019.

\bibitem{agrawal2018rewriting}
A.~Agrawal, R.~Verschueren, S.~Diamond, and S.~Boyd.
\newblock A rewriting system for convex optimization problems.
\newblock {\em Journal of Control and Decision}, 5(1):42--60, 2018.

\bibitem{arpon2020admm}
S.~Arp{\'o}n, T.~Homem-de Mello, and B.~Pagnoncelli.
\newblock An {ADMM} algorithm for two-stage stochastic programming problems.
\newblock {\em Annals of Operations Research}, 286(1):559--582, 2020.

\bibitem{banjac2017embedded}
G.~Banjac, B.~Stellato, N.~Moehle, P.~Goulart, A.~Bemporad, and S.~Boyd.
\newblock Embedded code generation using the {OSQP} solver.
\newblock In {\em 2017 IEEE 56th Annual Conference on Decision and Control
  (CDC)}, pages 1906--1911. IEEE, 2017.

\bibitem{bemporad1999robust}
A.~Bemporad and M.~Morari.
\newblock Robust model predictive control: A survey.
\newblock In {\em Robustness in identification and control}, pages 207--226.
  Springer, 1999.

\bibitem{bertsekas2011survey}
D.~Bertsekas.
\newblock Incremental gradient, subgradient, and proximal methods for convex
  optimization: A survey.
\newblock {\em Optimization for Machine Learning}, 2010(1-38):3, 2011.

\bibitem{bertsekas2011incremental}
D.~Bertsekas.
\newblock Incremental proximal methods for large scale convex optimization.
\newblock {\em Mathematical programming}, 129(2):163--195, 2011.

\bibitem{birge1988multicut}
J.~Birge and F.~Louveaux.
\newblock A multicut algorithm for two-stage stochastic linear programs.
\newblock {\em European Journal of Operational Research}, 34(3):384--392, 1988.

\bibitem{BoydVand04}
S.~Boyd and L.~Vandenberghe.
\newblock {\em Convex Optimization}.
\newblock Cambridge University Press, 2004.

\bibitem{caflisch1998monte}
R.~Caflisch.
\newblock Monte {C}arlo and quasi-{M}onte {C}arlo methods.
\newblock {\em Acta numerica}, 7:1--49, 1998.

\bibitem{camacho2013model}
E.~Camacho and C.~Bordons.
\newblock {\em Model predictive control}.
\newblock Springer science \& business media, 2013.

\bibitem{campo1987robust}
P.~Campo and M.~Morari.
\newblock Robust model predictive control.
\newblock In {\em 1987 American control conference}, pages 1021--1026. IEEE,
  1987.

\bibitem{CARLI2020373}
R.~Carli, G.~Cavone, N.~Epicoco, P.~Scarabaggio, and M.~Dotoli.
\newblock Model predictive control to mitigate the {COVID}-19 outbreak in a
  multi-region scenario.
\newblock {\em Annual Reviews in Control}, 50:373--393, 2020.

\bibitem{cutler1980dynamic}
C.~Cutler and B.~Ramaker.
\newblock Dynamic matrix control: A computer control algorithm.
\newblock In {\em Joint Automatic Control Conference}, volume~17, page~72,
  1980.

\bibitem{de2021risk}
W.~de~Oliveira.
\newblock Risk-averse stochastic programming and distributionally robust
  optimization via operator splitting.
\newblock {\em Set-Valued and Variational Analysis}, pages 1--31, 2021.

\bibitem{cvxpy_paper}
S.~Diamond and S.~Boyd.
\newblock {CVXPY}: A {P}ython-embedded modeling language for convex
  optimization.
\newblock {\em Journal of Machine Learning Research}, 2016.

\bibitem{ding2018104}
Y.~Ding, L.~Wang, Y.~Li, and D.~Li.
\newblock Model predictive control and its application in agriculture: A
  review.
\newblock {\em Computers and Electronics in Agriculture}, 151:104--117, 2018.

\bibitem{Domahidi2013ecos}
A.~Domahidi, E.~Chu, and S.~Boyd.
\newblock {ECOS}: {A}n {SOCP} solver for embedded systems.
\newblock In {\em European Control Conference (ECC)}, pages 3071--3076, 2013.

\bibitem{drgovna2020all}
J.~Drgo{\v{n}}a, J.~Arroyo, I.~Figueroa, D.~Blum, K.~Arendt, D.~Kim,
  E.~Oll{\'e}, J.~Oravec, M.~Wetter, D.~Vrabie, and L.~Helsen.
\newblock All you need to know about model predictive control for buildings.
\newblock {\em Annual Reviews in Control}, 2020.

\bibitem{engell2009online}
S.~Engell.
\newblock Online optimizing control: The link between plant economics and
  process control.
\newblock In {\em 10th International Symposium on Process Systems Engineering},
  volume~27, pages 79--86. Elsevier, 2009.

\bibitem{eren2017model}
U.~Eren, A.~Prach, B.~Ko{\c{c}}er, S.~Rakovi{\'c}, E.~Kayacan, and
  B.~A{\c{c}}{\i}kme{\c{s}}e.
\newblock Model predictive control in aerospace systems: {C}urrent state and
  opportunities.
\newblock {\em Journal of Guidance, Control, and Dynamics}, 40(7):1541--1566,
  2017.

\bibitem{felez2019model}
J.~Felez, Y.~Kim, and F.~Borrelli.
\newblock A model predictive control approach for virtual coupling in railways.
\newblock {\em IEEE Transactions on Intelligent Transportation Systems},
  20(7):2728--2739, 2019.

\bibitem{garcia1989model}
C.~Garcia, D.~Prett, and M.~Morari.
\newblock Model predictive control: Theory and practice—a survey.
\newblock {\em Automatica}, 25(3):335--348, 1989.

\bibitem{glynn1989importance}
P.~Glynn and D.~Iglehart.
\newblock Importance sampling for stochastic simulations.
\newblock {\em Management science}, 35(11):1367--1392, 1989.

\bibitem{golub1995stochastic}
B.~Golub, M.~Holmer, R.~McKendall, L.~Pohlman, and S.~Zenios.
\newblock A stochastic programming model for money management.
\newblock {\em European Journal of Operational Research}, 85(2):282--296, 1995.

\bibitem{grune2011nonlinear}
L.~Gr{\"u}ne and J.~Pannek.
\newblock Nonlinear model predictive control, 2011.

\bibitem{heirung2018stochastic}
T.~Heirung, J.~Paulson, J.~O’Leary, and A.~Mesbah.
\newblock Stochastic model predictive control -- how does it work?
\newblock {\em Computers \& Chemical Engineering}, 114:158--170, 2018.

\bibitem{holkar2010overview}
K.~Holkar and L.~Waghmare.
\newblock An overview of model predictive control.
\newblock {\em International Journal of Control and Automation}, 3(4):47--63,
  2010.

\bibitem{hovgaard2015model}
T.~Hovgaard, S.~Boyd, and J.~J{\o}rgensen.
\newblock Model predictive control for wind power gradients.
\newblock {\em Wind Energy}, 18(6):991--1006, 2015.

\bibitem{hu2021}
J.~Hu, Y.~Shan, J.~Guerrero, A.~Ioinovici, K.~Chan, and J.~Rodriguez.
\newblock Model predictive control of microgrids – an overview.
\newblock {\em Renewable and Sustainable Energy Reviews}, 136:110422, 2021.

\bibitem{kochenderfer2022algorithms}
M.~Kochenderfer, T.~Wheeler, and K.~Wray.
\newblock {\em Algorithms for decision making}.
\newblock Mit Press, 2022.

\bibitem{lazic2018data}
N.~Lazic, T.~Lu, C.~Boutilier, M.~Ryu, E.~J. Wong, B.~Roy, and G.~Imwalle.
\newblock Data center cooling using model-predictive control.
\newblock In {\em Proceedings of the Thirty-second Conference on Neural
  Information Processing Systems (NeurIPS-18)}, pages 3818--3827, Montreal, QC,
  2018.

\bibitem{lee2011model}
J.~Lee.
\newblock Model predictive control: Review of the three decades of development.
\newblock {\em International Journal of Control, Automation and Systems},
  9(3):415--424, 2011.

\bibitem{Lemaire1989}
B.~Lemaire.
\newblock The proximal algorithm.
\newblock {\em International Series of Numerical Mathematics}, pages 73 -- 87,
  1989.

\bibitem{linderoth2001implementing}
J.~Linderoth and S.~Wright.
\newblock Implementing decomposition algorithms for stochastic programming on a
  computational grid.
\newblock {\em Technical Report ANL/MCS-P909--0101}, 2001.

\bibitem{lucia2012new}
S.~Lucia, T.~Finkler, D.~Basak, and S.~Engell.
\newblock A new robust {NMPC} scheme and its application to a semi-batch
  reactor example.
\newblock {\em IFAC Proceedings Volumes}, 45(15):69--74, 2012.

\bibitem{maciejowski2002predictive}
J.~Maciejowski.
\newblock {\em Predictive control: with constraints}.
\newblock Pearson education, 2002.

\bibitem{maiworm2015scenario}
M.~Maiworm, T.~B{\"a}thge, and R.~Findeisen.
\newblock Scenario-based model predictive control: Recursive feasibility and
  stability.
\newblock {\em IFAC-PapersOnLine}, 48(8):50--56, 2015.

\bibitem{mattingley2012cvxgen}
J.~Mattingley and S.~Boyd.
\newblock {CVXGEN}: {A} code generator for embedded convex optimization.
\newblock {\em Optimization and Engineering}, 13(1):1--27, 2012.

\bibitem{mayne2014model}
D.~Mayne.
\newblock Model predictive control: Recent developments and future promise.
\newblock {\em Automatica}, 50(12):2967--2986, 2014.

\bibitem{mayne2005robust}
D.~Mayne, M.~Seron, and S.~Rakovi{\'c}.
\newblock Robust model predictive control of constrained linear systems with
  bounded disturbances.
\newblock {\em Automatica}, 41(2):219--224, 2005.

\bibitem{mesbah2016stochastic}
A.~Mesbah.
\newblock Stochastic model predictive control: {A}n overview and perspectives
  for future research.
\newblock {\em IEEE Control Systems Magazine}, 36(6):30--44, 2016.

\bibitem{moehle2019dynamic}
N.~Moehle, E.~Busseti, S.~Boyd, and M.~Wytock.
\newblock Dynamic energy management.
\newblock In {\em Large Scale Optimization in Supply Chains and Smart
  Manufacturing}, pages 69--126. Springer, 2019.

\bibitem{morari1999model}
M.~Morari and J.~Lee.
\newblock Model predictive control: past, present and future.
\newblock {\em Computers \& Chemical Engineering}, 23(4-5):667--682, 1999.

\bibitem{de2005stochastic}
D.~Munoz de~la Penad, A.~Bemporad, and T.~Alamo.
\newblock Stochastic programming applied to model predictive control.
\newblock In {\em Proceedings of the 44th IEEE Conference on Decision and
  Control}, pages 1361--1366. IEEE, 2005.

\bibitem{parikh2014proximal}
N.~Parikh and S.~Boyd.
\newblock Proximal algorithms.
\newblock {\em Foundations and Trends in optimization}, 1(3):127--239, 2014.

\bibitem{peni2020nonlinear}
T.~P{\'e}ni, B.~Csutak, G.~Szederk{\'e}nyi, and G.~R{\"o}st.
\newblock Nonlinear model predictive control with logic constraints for
  {COVID}-19 management.
\newblock {\em Nonlinear Dynamics}, 102(4):1965--1986, 2020.

\bibitem{PJM}
{PJM} data miner 2.
\newblock Available at \url{http://dataminer2.pjm.com/list}.

\bibitem{prekopa1995two}
A.~Pr{\'e}kopa.
\newblock Two-stage stochastic programming problems.
\newblock In {\em Stochastic Programming}, pages 373--423. Springer, 1995.

\bibitem{primbs2009dynamic}
J.~Primbs.
\newblock Dynamic hedging of basket options under proportional transaction
  costs using receding horizon control.
\newblock {\em International Journal of Control}, 82(10):1841--1855, 2009.

\bibitem{raimondo2009min}
D.~Raimondo, D.~Limon, M.~Lazar, L.~Magni, and E.~Camacho.
\newblock Min-max model predictive control of nonlinear systems: A unifying
  overview on stability.
\newblock {\em European Journal of Control}, 15(1):5--21, 2009.

\bibitem{rakovic2018handbook}
S.~Rakovi{\'c} and W.~Levine.
\newblock {\em Handbook of model predictive control}.
\newblock Springer, 2018.

\bibitem{rawlings2017model}
J.~Rawlings, D.~Mayne, and M.~Diehl.
\newblock {\em Model predictive control: {T}heory, computation, and design},
  volume~2.
\newblock Nob Hill Publishing Madison, WI, 2017.

\bibitem{richalet1978model}
J.~Richalet, A.~Rault, J.~L. Testud, and J.~Papon.
\newblock Model predictive heuristic control.
\newblock {\em Automatica (journal of IFAC)}, 14(5):413--428, 1978.

\bibitem{rockafellar1976monotone}
R.~Rockafellar.
\newblock Monotone operators and the proximal point algorithm.
\newblock {\em SIAM journal on control and optimization}, 14(5):877--898, 1976.

\bibitem{rockafellar2019progressive}
R.~Rockafellar.
\newblock Progressive decoupling of linkages in optimization and variational
  inequalities with elicitable convexity or monotonicity.
\newblock {\em Set-Valued and Variational Analysis}, 27(4):863--893, 2019.

\bibitem{rockafellar1991scenarios}
R.~Rockafellar and R.~Wets.
\newblock Scenarios and policy aggregation in optimization under uncertainty.
\newblock {\em Mathematics of operations research}, 16(1):119--147, 1991.

\bibitem{ruszczynski1986regularized}
A.~Ruszczy{\'n}ski.
\newblock A regularized decomposition method for minimizing a sum of polyhedral
  functions.
\newblock {\em Mathematical programming}, 35(3):309--333, 1986.

\bibitem{ruszczynski2003decomposition}
A.~Ruszczy{\'n}ski.
\newblock Decomposition methods.
\newblock {\em Handbooks in operations research and management science},
  10:141--211, 2003.

\bibitem{sand2004modeling}
G.~Sand and S.~Engell.
\newblock Modeling and solving real-time scheduling problems by stochastic
  integer programming.
\newblock {\em Computers \& chemical engineering}, 28(6-7):1087--1103, 2004.

\bibitem{sen1994network}
S.~Sen, R.~Doverspike, and S.~Cosares.
\newblock Network planning with random demand.
\newblock {\em Telecommunication systems}, 3(1):11--30, 1994.

\bibitem{sen1994solution}
S.~Sen, J.~Mai, and J.~L. Higle.
\newblock Solution of large scale stochastic programs with stochastic
  decomposition algorithms.
\newblock In {\em Large Scale Optimization}, pages 388--410. Springer, 1994.

\bibitem{shapiro2003monte}
A.~Shapiro.
\newblock Monte {C}arlo sampling approach to stochastic programming.
\newblock In {\em ESAIM: Proceedings}, volume~13, pages 65--73. EDP Sciences,
  2003.

\bibitem{stellato2020osqp}
B.~Stellato, G.~Banjac, P.~Goulart, A.~Bemporad, and S.~Boyd.
\newblock {OSQP}: {A}n operator splitting solver for quadratic programs.
\newblock {\em Mathematical Programming Computation}, 12(4):637--672, 2020.

\bibitem{topaloglou2008dynamic}
N.~Topaloglou, H.~Vladimirou, and S.~Zenios.
\newblock A dynamic stochastic programming model for international portfolio
  management.
\newblock {\em European Journal of Operational Research}, 185(3):1501--1524,
  2008.

\bibitem{tuck2021distributed}
J.~Tuck, S.~Barratt, and S.~Boyd.
\newblock A distributed method for fitting {L}aplacian regularized stratified
  models.
\newblock {\em Journal of Machine Learning Research}, 22(60):1--37, 2021.

\bibitem{tuck2021eigen}
J.~Tuck and S.~Boyd.
\newblock Eigen-stratified models.
\newblock {\em Optimization and Engineering}, pages 1--23, 2021.

\bibitem{tuck2021fitting}
J.~Tuck and S.~Boyd.
\newblock Fitting {L}aplacian regularized stratified {G}aussian models.
\newblock {\em Optimization and Engineering}, pages 1--21, 2021.

\bibitem{van1969shaped}
R.~Van~Slyke and R.~Wets.
\newblock L-shaped linear programs with applications to optimal control and
  stochastic programming.
\newblock {\em SIAM journal on applied mathematics}, 17(4):638--663, 1969.

\bibitem{wang2009fast}
Y.~Wang and S.~Boyd.
\newblock Fast model predictive control using online optimization.
\newblock {\em IEEE Transactions on control systems technology},
  18(2):267--278, 2009.

\bibitem{yu2013model}
Y.~Xi, D.~Li, and S.~Lin.
\newblock Model predictive control—status and challenges.
\newblock {\em Acta Automatica Sinica}, 39(3):222--236, 2013.

\bibitem{zenios1998dynamic}
S.~Zenios, M.~Holmer, R.~McKendall, and C.~Vassiadou-Zeniou.
\newblock Dynamic models for fixed-income portfolio management under
  uncertainty.
\newblock {\em Journal of Economic Dynamics and Control}, 22(10):1517--1541,
  1998.

\end{thebibliography}
\bibliographystyle{plain}
\end{document}